\documentclass[12pt, reqno, a4paper]{amsart}
\usepackage[margin=1.2in]{geometry}
\numberwithin{equation}{section}
\usepackage{amssymb,amsfonts,amsthm}
\usepackage[utf8]{inputenc}
\usepackage{listings}
\usepackage{bm}
\usepackage[hyperpageref]{backref}
\usepackage{esint}
\usepackage{bigints}
\usepackage{amsthm}
\usepackage{enumitem}
\usepackage{color}
\usepackage[bookmarks]{hyperref}
\usepackage{siunitx}
\usepackage{cite}
\usepackage[T1]{fontenc}
\theoremstyle{plain}
\usepackage{tikz}
\usetikzlibrary{shapes}
\usepackage{caption} 

\usepackage{mathtools}

\allowdisplaybreaks[4] 
\usepackage[T1]{fontenc}
\usepackage{amsthm}
\usepackage{enumitem}

\newtheoremstyle{myremark}{10pt}{10pt}{}{}{\bfseries}{.}{.5em}{}

\newtheorem{theorem}{Theorem}[section]
\newtheorem{lemma}[theorem]{Lemma}   
\newtheorem{prop}[theorem]{Proposition}

\newtheorem{definition}{Definition}[section]

\usepackage{hyperref}

\allowdisplaybreaks[4]

\begin{document}

\title[ FEM for a Stochastic Pseudo-Parabolic Equation ]{Strong convergence of finite element approximations for a fourth-order stochastic pseudo-parabolic equation with additive noise}
\address{Department of Mathematics and Statistics, Indian Institute of Technology Kanpur, Kanpur-208016, India}

\author{Suprio Bhar}
\email{Suprio Bhar: suprio@iitk.ac.in}

\author{Mrinmay Biswas }
\email{Mrinmay Biswas: mbiswas@iitk.ac.in}

\author{Mangala Prasad}
\email{Mangala Prasad: mangalap21@iitk.ac.in}
 
	

\date{}

\dedicatory{}

\begin{abstract}
In this article, we analyze semi-discrete finite element approximation and full discretization of a fourth-order stochastic
pseudo-parabolic equation in a bounded convex polygonal domain driven by additive Wiener noise. We use the finite element method for spatial discretization and the semi-implicit method for temporal discretization, and obtain strong convergence rates with respect to both the spatial and temporal mesh sizes. Numerical experiments are presented to support the theoretical convergence rates.
	\end{abstract}

\maketitle

\section{Introduction}\label{s1}
Pseudo-parabolic equations describe a variety of physical and engineering phenomena, including heat conduction in materials with memory, flow in porous media, and dispersive wave propagation \cite{BARENBLATT19601286,hassanizadeh1993thermodynamic,Juanes_2009,Aslan_2008,BBM1}. Equations involving mixed spatial and temporal derivatives exhibit behavior that differs from classical parabolic models and are commonly referred to as pseudo-parabolic. This class of equations was designated as “pseudo-parabolic” by Showalter and Ting in 1970 \cite{TW_Ting_1970}. This paper focuses on a fourth-order stochastic pseudo-parabolic equation with additive noise  
    \begin{equation}\label{eqn1.1}
        du-d(\Delta u)-(\Delta u-\Delta ^2u)dt =f(u)dt+dW,
    \end{equation}
    where $f$ is a nonlinear source term  and $\{W(t)\}_{t \geq 0}$ is an $L_2(\mathcal{O})$-valued Wiener process with covariance operator $Q$ on a filtered probability space $(\Omega, \mathcal{F},\{\mathcal{F}_t\}_{t \geq 0}, \mathbb{P})$ with respect to the filtration $\{\mathcal{F}_t\}_{ t \geq 0}$.

For the fourth-order nonlinear pseudo-parabolic equation \eqref{eqn1.1} without a noise term, global existence and asymptotic behavior of solutions have been extensively investigated \cite{Karlen_I_2009,Zhao1,Huafei_Di_2015,Huafei_2016} and the references therein. In addition to qualitative analysis, numerical approximations of pseudo-parabolic equations have also attracted considerable attention \cite{Keyan_2024,Amirali_2024,AmiyaK_Pani_2018,Yunxia_2025,AmiyaK_Pani_2019}. Various numerical methods, including finite difference schemes, spatial finite element discretizations, and backward Euler time-stepping methods, have been employed to analyze and approximate such models.

In many realistic applications, however, deterministic models are often influenced by random perturbations arising from environmental fluctuations, measurement errors, or unresolved physical mechanisms. This motivates incorporating stochastic effects into pseudo-parabolic equations via Wiener processes, yielding stochastic pseudo-parabolic models. Theoretical investigations concerning the existence, uniqueness, and asymptotic behavior of solutions to stochastic pseudo-parabolic equations have been reported in \cite{sotelo_2025,Sotelo_2024,cao2026parameter}. Moreover, stochastic fractional pseudo-parabolic equations have also been studied, and results on the existence and uniqueness of solutions can be found in \cite{Tuan_2023,Thach_2022}.

Despite the growing body of literature on the qualitative analysis of stochastic pseudo-parabolic equations, numerical approximation for such stochastic models remains relatively underdeveloped. In particular, to the best of our knowledge, strong convergence analysis for fully discrete finite element schemes combined with semi-implicit time discretization applied to fourth-order stochastic pseudo-parabolic equations has not yet been established. The presence of mixed spatial and temporal derivatives, together with stochastic forcing, introduces additional analytical challenges compared with classical stochastic parabolic equations.

By contrast, finite element approximations for stochastic heat, wave, and Schr\"odinger equations have been extensively investigated in the literature; see, for example, \cite{Yan1,Yan2,kovacs1,kovacs2,chuchuhong,finite}. Moreover, semi-implicit time discretization schemes have recently attracted considerable attention due to their stability and efficiency; see, e.g., \cite{yan3,time16,Tretyakov}. These developments provide a solid foundation for constructing and analyzing fully discrete finite element schemes for the fourth-order stochastic pseudo-parabolic equation considered in this work, which involves the interplay of fourth-order spatial derivatives, nonlinear effects, and stochastic perturbations.

In this paper, we develop a finite element approximation and a semi-implicit time discretization for the fourth-order stochastic pseudo-parabolic equation with respect to both spatial and temporal variables, derive the convergence rates, and investigate the well-posedness of its solutions. Due to the involvement of the dissipative term $\Delta ^2u$ and mixed derivative, finding the rate of convergence and proving the existence of the solution of the equation \eqref{eqn1.1} is quite challenging. To address this problem, we transform the fourth-order stochastic pseudo-parabolic equation into a stochastic coupled parabolic–elliptic system by introducing a new variable $v=u-\Delta u$, i.e
 \begin{equation}\label{eqn1.2}
     \begin{cases}
        dv-\Delta vdt=f(u)dt+dW,\\
        v=u-\Delta u.
        \end{cases}
\end{equation}
\textbf{Main contributions.}
The main contributions of this paper are summarized as follows.
\begin{itemize}
    \item We establish the existence and uniqueness of mild solutions for the fourth-order stochastic pseudo-parabolic equation driven by additive Wiener noise and derive suitable regularity estimates that serve as the foundation for the numerical analysis.
    \item We analyze both the semi-discrete finite element approximation in space and the fully discrete scheme that combines a finite element spatial discretization with a semi-implicit time-stepping method. Strong convergence rates are derived with respect to the spatial and temporal mesh sizes.
    \item Numerical experiments with explicit examples are presented to validate the theoretical convergence results
\end{itemize}
This paper is organized as follows. 
In Section~\ref{S2}, we introduce the problem setting and state the main results. 
Section~\ref{s3} is devoted to preliminaries, including basic definitions, regularity properties, finite element approximations, and results concerning semigroup approximations. 
In Section~\ref{s4}, we present the proofs of the results stated in Section~\ref{S2}. 
Finally, in Section~\ref{S5}, we provide a numerical example that validates the theoretical error bounds.

Throughout the paper, the constant $C$ denotes a generic positive constant independent of the discretization parameters. Its value may vary from line to line.
\section{Setup and Main Results} \label{S2}
Let $\mathcal{O} \subset \mathbb{R}^d, \, d=1,2,3$ be a bounded convex polygonal domain with boundary $\partial\mathcal{O}$. We consider the fourth-order stochastic pseudo-parabolic  equation 
\begin{equation}\label{eqn1.3}
    \begin{cases}
        du-d(\Delta u)-(  \Delta u-\Delta ^2u)dt=f(u)dt+dW \quad \text{ in } (0,\infty) \times \mathcal{O},\\
        u=\Delta u=0 \hspace{5.3cm}\qquad \text{ on } (0,\infty) \times \partial\mathcal{O},\\
        u(0,x)=u_0(x) \hspace{5.7cm} \text{ in } \mathcal{O},
    \end{cases}
\end{equation}
 where $f:L_2(\mathcal{O})\to L_2(\mathcal{O})$ is a globally Lipschitz function and  $u_0$ is an  $\mathcal{F}_0$ measurable random variable. To analyze the above system, we introduce the change of variable $v=u-\Delta u$, which converts equation \eqref{eqn1.3} into a stochastic coupled   parabolic-elliptic  system
\begin{equation}\label{eqn1.4}
    \begin{cases}
        dv-\Delta vdt=f(u)dt+dW \quad \text{ in }(0,\infty)\times \mathcal{O},\\
        v=u-\Delta u \hspace{3.1cm} \text { in } (0,\infty)\times \mathcal{O},\\
         v=u=0  \hspace{3.3cm}\text{ on }(0,\infty) \times\partial\mathcal{O},\\
         v(0,x)=v_0(x ) \hspace{2.5cm}\text{ in } \mathcal{O},
          \end{cases}
\end{equation}
where $v_0=u_0-\Delta u_{0}$. We observe that in the above system, the equation in $v$ is a parabolic equation coupled with a nonlinear forcing term $f(u)$ and a stochastic perturbation. The second equation in $u$ is an elliptic equation coupled with a linear forcing term $v$. In the elliptic equation, $u(t)$ plays the role of a parameter. Thus, it is a family of elliptic equations. We note that solving for $v$ with the initial condition 
$v_0$ yields the solution 
$u$ associated with the initial condition 
$u_0$. Our analysis will focus on 
$v$ to derive the solution and establish the convergence rate of the finite element approximation as well as semi-implicit time discretization for the equation \eqref{eqn1.3}.

 First, we consider the stochastic parabolic component of the equation \eqref{eqn1.4}
    \begin{equation}\label{eqn1.5}
       dv(t)+Av(t)dt=f(u(t))dt+dW(t), \quad t >0; \, v(0)=v_0,
    \end{equation}
    where $A=-\Delta$  denotes the Laplace operator with domain $D(A)=H^2(\mathcal{O}) \cap H^1_0(\mathcal{O}) $. The operator $-A$ generates the $C_0$ semigroup $\{E(t)\}_{t \geq 0}$.\\
    We say that a process $\{v_t\}_{t \geq 0}$ is a mild solution of \eqref{eqn1.5}, if it is $\mathcal{F}_t$ adapted and satisfies the integral equation $\mathbb{P}$ a.s.
    \begin{equation}\label{eqn1.6}
        v(t) =E(t)v_0+\int_0^t E(t-\tau)f(u(\tau))d\tau+\int_0^tE(t-\tau)dW(\tau),\quad \text{for } t \geq 0.
    \end{equation}
  Let $V_h$ be a finite-dimensional subspace of $H_0^1(\mathcal{O})$. It is defined in Subsection \ref{s2.2}. The finite element approximation of equation \eqref{eqn1.6} is to find $v_h \in V_h$ such that 
\begin{equation}\label{eqn1.7}
dv_h(t)+A_hv_h(t)dt=\mathcal{P}_hf(u_h(t))dt+\mathcal{P}_hdW(t) ,\quad  t>0, \, v_h(0)=\mathcal{P}_hv_0=v_{0,h},
    \end{equation}
     where $A_h$ is the discrete Laplace operator. The discrete operator $-A_h$ generates the continuous semigroup $\{E_h(t)\}_{t \geq 0}$. Operator  $\mathcal{P}_h$ denotes the $L_2$ projection onto the finite element space $V_h$, and $u_h(t)$ is a finite element approximation of equation \eqref{eqn1.3}. Using the transformation  $v=u+Au=(I+A)u$,  the function $u_h(t)$ can be expressed in terms of $v_h(t)$ as follow 
    \begin{equation} \label{eqn1.8}
     u_h(t)=(\mathcal{P}_h +A_h)^{-1}v_h(t) \quad \text{ for } t \geq 0.
  \end{equation}
The discrete operator $ A_h$ and the projection $\mathcal{P}_h$ are defined in  Subsection \ref{s2.2}.
In analogy with \eqref{eqn1.6}, the mild solution of the discretized equation \eqref{eqn1.7} is expressed as $\mathbb{P}$ a.s. 
     \begin{equation}\label{eqn1.9}
         v_h(t)=E_h(t)v_{0,h}+\int_0^tE_h(t-\tau)\mathcal{P}_hf(u_h(\tau))d\tau+ \int_0^t E_h(t-\tau)\mathcal{P}_hdW(\tau), \quad \text{ for } t \geq 0.
     \end{equation}
     First, we discuss the existence and uniqueness of a mild solution of the equation \eqref{eqn1.3}.
\begin{prop}\label{prop1}
Let $f$ be a globally Lipschitz function, $Q$ be trace class operator satisfying  $\|A^{(\beta-1)/2}Q
^{1/2}\|_{HS}< \infty$ for $\beta \geq 0$, and $u_0$ is an $\mathcal{F}_0$ measurable $H^2(\mathcal{O})$-valued random variable. Then there is a unique mild solution $\{u(t)\}_{t \geq 0}$ to the equation \eqref{eqn1.3}. In particular $\|u(t)\|_{L_2(\Omega,L_2(\mathcal{O}))} \leq C$ for $0 \leq t \leq T$, where $C=C(T,f)$.
\end{prop}
Before determining the convergence rate of the finite element approximation, it is necessary to establish an a priori estimate for  $v_h(t)$.
\begin{lemma} \label{lemma1}
Let \( Q \) be a trace-class operator, \( f \) a globally Lipschitz continuous function, and let \( v_0 \in L_2(\Omega, L_2(\mathcal{O})) \). Then, for   \( 0 \leq t \leq T \), there exists a constant \( C = C(T, f, Q) \) such that
\[
\|v_h(t)\|_{L_2(\Omega, L_2(\mathcal{O}))} \leq C.
\]
\end{lemma}
By applying Lemma \ref{lemma1} in combination with finite element techniques developed for deterministic parabolic problems (see, for instance, \cite[Theorem 3.5]{thomee} ), as recalled in Subsection \ref{s2.2}, we derive the following convergence rate for the fourth-order stochastic pseudo-parabolic equation \eqref{eqn1.3}. In this analysis, we employ the fractional Sobolev space $ \dot{H}^{\beta}$, defined in Section \ref{s3}.
\begin{theorem}\label{thm1}
Let \( \{u(t)\}_{t \geq 0} \) denote the mild solution of equation \eqref{eqn1.3}, and let \( \{u_h(t)\}_{t \geq 0} \) be its finite element approximation given by equation \eqref{eqn1.8}. Suppose that \( \|A^{(\beta - 1)/2} Q^{1/2}\|_{\mathrm{HS}} < \infty \) for some \( \beta \in [0, 2) \), and  $v_0 \in L_2(\Omega,\dot{H}^{\beta})$. Then there exists a constant \( C = C(v_0, T, Q, f) > 0 \) such that for  \(0 \leq \beta < 2 \),
\[
\|u(t) - u_h(t)\|_{L_2(\Omega, L_2(\mathcal{O}))} \leq C h^{\beta}, \quad 0\leq t\leq T.
\]
\end{theorem}
Based on the semi-discrete finite element approximation discussed above, we introduce an additional temporal discretization. This leads to a fully discrete formulation, for which we derive analogous stability and convergence results.

We assume that the function $f:L_2(\mathcal{O}) \to L_2(\mathcal{O})$ satisfies the following global Lipschitz and growth conditions:
\begin{enumerate}[label=(\roman*),ref=\roman*]
    \item $\|f(x_1)-f(x_2)\| \leq K_f\|x_1-x_2\|,\quad x_1,\,x_2\in L_2(\mathcal{O}).$ \label{cond1}
    \item $\|A^{(\beta-1)/2}f(x)\| \leq C\|x\|,\quad \forall x\in L_2(\mathcal{O})\text{ and for some } \beta \in [0,1].$ \label{cond2}
\end{enumerate}
To obtain full discretization, we need regularity in time of the solution \eqref{eqn1.3}.
\begin{lemma}\label{lemma2}
    Assume that the growth condition  \eqref{cond2} holds. Let $u$ be the mild solution of \eqref{eqn1.4}. Then for $0 \leq \gamma < \beta \leq1$, we have 
    \[
    \mathbb{E}\|u(t_2)-u(t_1)\|^2\leq C(t_2-t_1)^{\gamma}.
    \]
\end{lemma}
Let $k$ be a time step and $t_n=kn$. The fully discrete approximation is defined by the semi-implicit time-stepping scheme of the semi-discrete equation \eqref{eqn1.9}, i.e., find $V^n \in V_h$, such that 
\begin{equation}\label{eqn2.8}
   V^n-V^{n-1}+kA_hV^n= k\mathcal{P}_hf(U^{n-1})+\int _{t_{n-1}}^{t_n}\mathcal{P}_hdW(\tau), \, n\geq 1,\, V^0=\mathcal{P}_hv_0,
\end{equation}
where $U^n$ is full discretization of equation \eqref{eqn1.3} and it is expressed in terms of $V^n$ as follow
\begin{equation}\label{fullydu}
    U^n=(\mathcal{P}_h+A_h)^{-1}V^n, \text{ for } \,n \geq 0.
\end{equation}
 The fully discrete scheme \eqref{eqn2.8} satisfies the mean-square stability property, as stated in the following proposition.
 \begin{prop} \label{prop2}
     The numerical solution of the fully discrete equation \eqref{eqn2.8} is mean square stable, i.e.
     \[
     \mathbb{E}\|V^n\|^2\leq C(T,Q,v_0,f), \quad n=1,2,\cdots.
     \]
 \end{prop}
In the next result, we state our second main result, which is related to the full discretization. With $R(\lambda)=1/(1+\lambda)$, we can rewrite \eqref{eqn2.8} in the form 
\begin{equation}\label{eqn2.9}
    \begin{split}
        &V^{n}=R(kA_h)V^{n-1}+R(kA_h)k\mathcal{P}_hf(U^{n-1})+\int_{t_{n-1}}^{t_n}R(kA_h)\mathcal{P}_hdW(\tau),\quad n \geq 1,\,\\
        &V^0=\mathcal{P}_hv_0.
    \end{split}
\end{equation}
By combining the temporal continuity of the solution with its mean-square stability, we derive the mean-square error estimate for the fully discrete approximation of the fourth-order stochastic pseudo-parabolic equation.
\begin{theorem}\label{thm2}
Let $u(t_n)$ denote the exact solution of \eqref{eqn1.3} and $U^n$ its fully discrete approximation defined by \eqref{fullydu}. 
Assume that $f$ satisfies conditions \eqref{cond1} and \eqref{cond2}. 
Suppose further that $v_0 \in L_2(\Omega,\dot{H}^{\beta})$ and $\|A^{(\beta-1)/2}Q^{1/2}\|_{HS} < \infty$ for $0 \leq \beta \leq 1$. 
Then there exists a constant $C=C(f,T,Q,v_0)>0$ such that, for $t_n \in [0,T]$ and $0 \leq \gamma < \beta \leq 1$, the following error estimate holds:
\[
\|U^n - u(t_n)\|_{L_2(\Omega,L_2(\mathcal{O}))} \leq C\big(k^{\gamma/2} + h^{\beta}\big).
\]
\end{theorem}

\section{ Preliminaries} \label{s3}
	Let $(U_1,(\cdot,\cdot)_{U_1})\text{ and }(U_2,(\cdot,\cdot)_{U_2})$ be separable Hilbert spaces with corresponding norms $\|\cdot\|_{U_1} \text{ and } \|\cdot\|_{U_2}$. We suppress the subscripts when the spaces in use are clear from the context. Let $\mathcal{L}(U_1,U_2)$ denote the space of bounded linear operators from $U_1 \text{ to } U_2$, and $\mathcal{L}_2(U_1,U_2)$ the space of Hilbert-Schmidt operators, endowed with norm $\|\cdot\|_{\mathcal{L}_2(U_1,U_2)}$. That is, $T \in \mathcal{L}_2(U_1,U_2) $ if $T \in \mathcal{L}(U_1,U_2)$ and 
	\[
	\|T\|^2_{\mathcal{L}_2(U_1,U_2)}:=\sum_{j=1}^{\infty}\|Te_j\|_{U_2}^2< \infty,
	\]
	where $\{e_j\}_{j=1}^{\infty}$ is an orthonormal basis of $U_1$. If $U_1=U_2$, we write $\mathcal{L}(U_1)=\mathcal{L}(U_1,U_1)$ and $HS= \mathcal{L}_2(U_1,U_1)$. It is a well-established fact that if
 $T_1 \in \mathcal{L}(U_1) \text{ and } T_2 \in \mathcal{L}_2(U_1,U_2)$, then their composition  $T_2T_1 \text{ also belongs to }\mathcal{L}_2(U_1,U_2)$. Moreover, the following norm bound holds:
	\[
	\|T_2T_1\|_{\mathcal{L}_2(U_1,U_2)} \leq \|T_2\|_{\mathcal{L}_2(U_1,U_2)}\|T_1\|_{\mathcal{L}(U_1)}.
	\]
	
	Let $\left(\Omega, \mathcal{F}, \{\mathcal{F}_t\}_{t \geq 0},\mathbb{P} \right)$ be a filtered probability satisfying the usual conditions. We denote by $L_2(\Omega, U_2)$  the space consisting of square-integrable random variables that take values in $U_2$. The norm associated with this space is given by
	\[
	\|v\|_{L_2(\Omega,U_2)}=\mathbb{E}(\|v\|^2_{U_2})^{1/2}=\left(\int_{\Omega}\|v(\omega)\|_{U_2}^2 dP(\omega)\right)^{1/2},
	\]
	where $\mathbb{E}$ is expectation. 
  \begin{definition}[\cite{rockner}]

   Consider an operator $Q \in \mathcal{L}(U_1)$ that is self-adjoint, positive, and semidefinite, satisfying $Tr(Q)< \infty$, where $Tr(Q)$ denotes its trace. A stochastic process $\{W(t)\}_{t \geq 0}$ taking values in $U_1$ is called  $Q$-Wiener process with respect to the filtration $\{\mathcal{F}_t\}_{t \geq 0}$ if it satisfies the following conditions:
	\begin{enumerate}
		\item $W(0)=0$  almost surely,
		\item $W$  has continuous sample paths with probability one,
		\item The increments of $W$ are independent,
		\item  For any $ 0\leq s \leq t$, the increment $W(t)-W(s), $  follows a Gaussian distribution in $U_1$ with mean zero and covariance operator
		 $(t-s)Q$.
	\end{enumerate}
\end{definition} 
Moreover the stochastic process  $\{W(t)\}_{t \geq 0}$ has the following representation
	\begin{equation}\label{eqn7}
		W(t)=\sum_{j=1}^{\infty} \gamma _j^{1/2} \beta_j(t)e_j,
	\end{equation}
	where $\{(\gamma_j,e_j)\}_{j=1}^{\infty}$ are the eigenpairs of $Q$, with $\{e_j\}$  forming an orthonormal basis, and
     $\{\beta _j\}_{j=1}^{\infty}$
	is a collection of independent standard Brownian motions in $\mathbb{R}.$ It can be verified that the series representation in \eqref{eqn7} converges in  $L_2(\Omega, U_1)$. Specifically, for any $t \geq 0$, we have
		\[
		\|W(t)\|_{L_2(\Omega,U_1)}^2= \mathbb{E}\left(\left\|\sum _{j=1} ^{\infty}\gamma_j ^{1/2} e_j \beta_j(t) \right\|_{U_1}^2\right)
	\]
Applying the linearity of expectation and orthonormality of $\{e_j\}$ this expression simplifies to
\[
    \sum_{j=1}^{\infty}\gamma_j\mathbb{E}(\beta_j(t))^2	 
    \]	
    Since each $\beta_j(t)$   is a standard Brownian motion, we use the property $\mathbb{E}\beta_j(t)^2=t$  to obtain
     \[
    \|W(t)\|_{L_2(\Omega,U_1)}^2= t \sum_{j=1}^{\infty}\gamma_j=t \text{Tr}(Q).
     \]   
        
 We consider here a specific instance of  It\^o's integral, where the integrand is a deterministic function.
	Let $\Phi : [0, \infty) \to  \mathcal{L}(U_1, U_2) $ be a strongly measurable function that satisfies the integrability condition
	\begin{equation}\label{eqn8}
		\int_0^t\|\Phi(\tau)Q^{1/2}\|^2_{\mathcal{L}_2(U_1, U_2) }d\tau < \infty.
	\end{equation}
	Under this assumption, the stochastic integral $\int_0^t\Phi(\tau)\,dW(\tau)$is well-defined, and It\^o's isometry isometry holds:
	\begin{equation} \label{eqn9}
		\left\|\int_0^t\Phi(\tau)\,dW(\tau)\right\|^2_{L_2(\Omega,U_2)}=\int_0^t \|\Phi (\tau) Q^{1/2}\|^2_{\mathcal{L}_2(U_1,U_2)}\,d\tau.
	\end{equation}
	
In a more general setting, consider a self-adjoint, positive, and semidefinite operator $Q \in \mathcal{L}(U_1)$  with eigenpairs $\{(\gamma_j, e_j )\}_{j=1}^{\infty}$. If $Q$ is not of trace class, meaning that  $Tr(Q) = \infty$,  then the series expansion given in \eqref{eqn7} fails to converge in $L_2(\Omega, U_1)$. However, it converges in a suitably chosen (usually larger) Hilbert space, and the stochastic integral $\int ^t _0 \Phi(\tau) dW(\tau) $ can still be defined, and the isometry \eqref{eqn9} holds, as long as \eqref{eqn8} is satisfied. In this case, $W$ is called a cylindrical Wiener process (\cite{prato}). In particular, we may have $Q = I$ (the identity operator).	
	\subsection{Abstract framework and regularity }\label{s2.1}
    Let us consider the Hilbert space $L_2(\mathcal{O})$ of square integrable Lebesgue measurable functions defined on $\mathcal{O}$ with usual inner product $(\cdot,\cdot)$ and norm $\|\cdot\|.$
	Let us define the Laplace operator $({A}, D({A}))$ in $L_2(\mathcal{O})$
    as
    \begin{align*}
        D({A}):&=H^2(\mathcal{O}) \cap H^1_0(\mathcal{O}),\\
        {A}u :& =-\Delta u,\quad u
        \in D({A}).
    \end{align*} 
    From spectral analysis (see for e.g. \cite{Yosida}), we know that $({A}, D({A}))$ has eigenpairs $\{(\lambda_j,\phi_j)\}_{j=1}^{\infty}$ satisfying the following:
    $$
    0<\lambda_1< \lambda_2 \leq \lambda_3\leq \cdots,\quad \lim_{j\to\infty} \lambda_j=\infty \quad\text{and}\quad \|\phi_j\|=1,\,\forall j\geq 1.
    $$ We introduce the following fractional Sobolev space using the above spectral decomposition of $({A}, D({A}))$
    as 
	\begin{equation*}
		\dot{H} ^{s}:=D({A} ^{s/2}),\hspace{1cm}\|v\|_{s}:=\|{A}^{s/2}v\|=\left(\sum_{j=1}^{\infty}\lambda^{s}_j(v,\phi_j)^2\right)^{1/2},\quad s \in \mathbb{R},\quad v\in \dot{H} ^{s}. 
	\end{equation*}
    We have $\dot{H} ^{s_2} \subset \dot{H} ^{s_1} \text{ for } s_1 \leq s_2.$ The spaces $ \dot{H} ^{-s}$ can be identified with the dual space $( \dot{H} ^{s})^{\star} \text{ for } s >0;$ see\cite{thomee}. We note that $\dot{H} ^{0}=L_2(\mathcal{O})=H,\text{ } \dot{H} ^{1}={H} ^{1}_0(\mathcal{O}), \text{ } \dot{H} ^{2}= D({A})=H^2(\mathcal{O}) \cap H^1_0(\mathcal{O})$ with equivalent norms . 
	\subsection{Finite element approximations}\label{s2.2}
	Let $\mathcal{T}_h$ represent a regular family of triangulations of 
 the domain $\mathcal{O}$. The parameter $h$ is defined as  $$h=\max_{K\in \mathcal{T}_h}h_K,$$
 where  $h_K=\text{diam}(K)\text{ and, } K \text{ is triangle in }\mathcal{O}$. Define $V_h$ as the space of continuous, piecewise linear functions associated with $\mathcal{T}_h$ that vanish on $\partial \mathcal{O}$.  Consequently, we have $V_h \subset H^1_0(\mathcal{O})=\dot{H}^1.$ 
	
	Since $\mathcal{O}$  is assumed to be convex and polygonal, the triangulations can be precisely aligned with $\partial \mathcal{O}$. Moreover, this ensures the elliptic regularity property $\|w\|_{H^2(\mathcal{O})} \leq C\|Aw\|$ for $w \in D({A})$, see \cite{pgrisvard}. We now recall fundamental results from the finite element theory (see \cite{ciarlet, Brenner}). The norms are denoted as $\|\cdot\|_s=\|\cdot\|_{\dot{H}^s}$.
	
	Consider the orthogonal projection operators  $\mathcal{P}_h:\dot{H}^0 \to V_h $ which is defined as
	\[
	(\mathcal{P}_hw,\psi)=(w, \psi).
	\]
	To introduce a discrete analogue of the norm 
$\|\cdot\|_{s}$, we define 
	\[
	\|w_h\|_{h,s}=\|A _h^{s/2}w\|, \hspace{.3cm} w_h\in V_h, s \in \mathbb{R},
	\]
	where the discrete Laplacian  $A_h: V_h \to V_h$  
  is given by
	\[
	(A_hw_h,\psi)=(\nabla w_h,\nabla \psi) \hspace{.2cm} \forall \psi \in V_h.
	\]
Next, we recall some standard results from deterministic finite element analysis and operator semigroup theory. 
The operator $-A$ generates an analytic semigroup $\{E(t)\}_{t \geq 0}$ on $H$, while its discrete counterpart $-A_h$, obtained by finite element discretization, generates the analytic semigroup $\{E_h(t)\}_{t \geq 0}$ on $V_h$. 
For the deterministic parabolic finite element problem, the following error estimate holds (see, e.g., \cite[Theorem~3.5]{thomee}):
\begin{equation}\label{destimate 1}
   \|(E_h(t)\mathcal{P}_h - E(t))\psi\| \leq C h^{\beta} t^{-\beta/2} \|\psi\|, 
   \quad 0 \leq \beta \leq 2, \; \psi \in H.
\end{equation}
In particular, from \eqref{destimate 1} we have 
\begin{equation}\label{destimate 2}
   \|E_h(t)\mathcal{P}_h\|\leq C,\quad t\geq 0, \, h>0. 
\end{equation}

\section{Proof of Main Results}\label{s4}
In this section, we prove the results that are stated in Section \ref{S2}.
\begin{proof}[\bf Proof of Proposition \ref{prop1}] We consider the stochastic coupled   parabolic-elliptic equation \eqref{eqn1.4}. From the change of variable $v(t)= (I+A)u(t)$ that is same as $u(t)=(I+A)^{-1}v(t)$ for $t \geq 0$ and equation \eqref{eqn1.6}, we obtain
\begin{equation}\label{eqn4.1}
   v(t) =E(t)v_0+\int_0^t E(t-\tau)f((I+A)^{-1}v(\tau))d\tau+\int_0^tE(t-\tau)dW(\tau),\quad \text{for } t \geq 0.
\end{equation}
The function $f$ is globally Lipschitz, and $(I+A)^{-1}$, being the resolvent of the Laplacian operator, is a bounded linear operator. Consequently, the mapping 
 $v(t) \mapsto f((I+A)^{-1}v(t))$ is also globally Lipschitz. Therefore, equation \eqref{eqn1.5} admits a unique mild solution (see \cite{kovacs1}, Lemma 3.1, and \cite{dprato1}, Theorem 2.3). This, in turn, yields the unique mild solution to equation \eqref{eqn1.3}, given by $u(t)=(I+A)^{-1}v(t)$ for $t \geq 0$. The boundedness of the solution $u(t)$ follows directly from the boundedness of both $v(t)$ and the resolvent operator $(I+A)^{-1}.$
\end{proof}
\begin{proof}[\bf Proof of Lemma \ref{lemma1}]
  We employ Gronwall’s inequality and It\^{o}'s isometry to estimate \( v_h(t) \). From equation \eqref{eqn1.9}, we obtain:
    \[
    \begin{split}
       \mathbb{E}\|v_h(t)\|^2  &=\mathbb{E}\left\|E_h(t)v_{0,h}+\int_0^tE_h(t-\tau)\mathcal{P}_hf(u_h(\tau))d\tau+ \int_0^t E_h(t-\tau)\mathcal{P}_hdW(\tau)\right\|^2\\
       &\leq 3\mathbb{E}\|E_h(t)v_{0,h}\|^2+3\mathbb{E}\left\|\int_0^tE_h(t-\tau)\mathcal{P}_hf(u_h(\tau))d\tau\right\|^2\\
       &+3\mathbb{E}\left\|\int_0^t E_h(t-\tau)\mathcal{P}_hdW(\tau)\right\|^2\\
       & \leq 3\mathbb{E}\|E_h(t)\mathcal{P}_hv_0\|^2+3t\mathbb{E}\int_0^t\left\|E_h(t-\tau)\mathcal{P}_hf(u_h(\tau))\right\|^2d\tau\\
       &+3\int_0^t \left\|E_h(t-\tau)\mathcal{P}_hQ^{1/2}\right\|_{HS}^2d\tau\\
       & \leq  3C\mathbb{E}\|v_0\|^2+3Ct\mathbb{E}\int _0^t\|f(u_h(\tau))\|^2d\tau+3Ct \text{Tr}(Q)\\
       & \leq (3C\mathbb{E}\|v_0\|^2+3Ct+3Ct\text{Tr}(Q))+3Ct\int _0^t\mathbb{E} \|u_h(\tau)\|^2d\tau.
    \end{split}
   \]
  Thus, we obtain:
   \begin{equation}\label{errorv}
       \mathbb{E}\|v_h(t)\|^2\leq B(t)+3Ct\int_0^t \mathbb{E}\|u_h(\tau)\|^2d\tau,
   \end{equation}
   where \(B(t)=3C\mathbb{E}\|v_0\|^2+3Ct+3Ct\text{Tr}(Q)\). Next, we establish the boundedness of \(\mathbb{E}\|u_h(t)\|^2 \) for  $t \geq 0$. From the equation \eqref{eqn1.8} and \eqref{errorv}, we have 
   \[
   \begin{split}
      \mathbb{E} \|u_h(t)\|^2&=\mathbb{E}\|(\mathcal{P}_h+A_h)^{-1}v_h(t)\|^2\\
      & \leq \mathbb{E}\|v_h(t)\|^2\\
      & \leq B(t)+3Ct\int_0^t \mathbb{E}\|u_h(\tau)\|^2d\tau.
   \end{split}
   \]
  Applying Gronwall’s inequality yields: 
   \[
     \mathbb{E} \|u_h(t)\|^2 \leq B(t)e^{3Ct^2}.
   \]
  Substituting this bound into inequality \eqref{errorv} gives:
   \[
   \mathbb{E}\|v_h(t)\|^2\leq B(t)+3Ct\int_0^t B(\tau) e^{3C\tau^2}d\tau.
   \]
   This completes the proof of Lemma \ref{lemma1}.
\end{proof}
Now we prove our first main result of the paper. First, we find the error between $v(t)$ and $v_h(t)$. Using this error, we obtain the rate of convergence of the finite element approximation of the fourth-order stochastic pseudo-parabolic equation \eqref{eqn1.3}.
\begin{proof}[\bf Proof of Theorem \ref{thm1}]
We use the representations \eqref{eqn1.6} and \eqref{eqn1.9}, the error between the exact and approximate solutions for 
     $t \geq 0$ is given by
         \[
         \begin{split}
             e_v(t):&=\|v(t)-v_h(t)\|_{L_2(\Omega, H)}\\
             &=\bigg\|E(t)v_0-E_h(t)v_{0,h}+\int_0^t(E(t-\tau)-E_h(t-\tau)\mathcal{P}_h)dW(\tau)\\
             & \qquad +\int_0^t\left(E(t-\tau)f(u(\tau))-E_h(t-\tau)\mathcal{P}_hf(u_h(\tau))\right)d\tau \bigg\|_{L_2(\Omega, H)}\\
             & \leq \left\|E(t)v_0-E_h(t)v_{0,h}+\int_0^t(E(t-\tau)-E_h(t-\tau)\mathcal{P}_h)dW(\tau)\right \|_{L_2(\Omega, H)}\\
             &+ \left \|\int_0^t\left(E(t-\tau)f(u(\tau))-E_h(t-\tau)\mathcal{P}_hf(u_h(\tau))\right)d\tau \right\|_{L_2(\Omega, H)}\\
             &=:I_1(t)+I_2(t).
         \end{split}
         \]
    By Proposition 3.2 of \cite{kovacs1}, we have the following estimate for $I_1(t)$
      \begin{equation}\label{eqn3.3}
          I_1(t)\leq Ch^{\beta }\left(\|v_0\|_{L_2(\Omega, \dot H^{\beta})}+\|A^{(\beta-1)/2}Q^{1/2}\|\right).
      \end{equation}
     For the $I_2(t)$, we use Lipschitz continuity property of the function $f$ and the deterministic finite element error estimate \eqref{destimate 1} and \eqref{destimate 2} to get 
     \[
     \begin{split}
         I_2(t)&= \left \|\int_0^t\left(E(t-\tau)f(u(\tau))-E_h(t-\tau)\mathcal{P}_hf(u_h(\tau))\right)d\tau \right\|_{L_2(\Omega, H)}\\
         &=\bigg\| \int _0^t \bigg[E(t-\tau)f(u(\tau))-E_h(t-\tau)\mathcal{P}_hf(u(\tau))\\
         & \hspace{2.0cm} + E_h(t-\tau)\mathcal{P}_hf(u(\tau))-E_h(t-\tau)\mathcal{P}_hf(u_h(\tau))\bigg]d\tau\bigg\|_{L_2(\Omega, H)}\\
 & \leq \int _0^t \left\|\left(E(t-\tau)-E_h(t-\tau)\mathcal{P}_h\right)f(u(\tau))\right\|_{L_2(\Omega, H)}d\tau\\
         &+\int_0^t \left\|E_h(t-\tau )\mathcal{P}_h(f(u(\tau))-f(u_h(\tau))) \right\|_{L_2(\Omega, H)}d\tau\\
         & \leq Ch^{\beta} \int _0^t (t-\tau)^{-\beta/2 }\|f(u(\tau))\|_{L_2(\Omega, H)}d\tau\\
         &+\int_0^t \|f(u(\tau))-f(u_h(\tau))\|_{L_2(\Omega, H)}d\tau\\
         \end{split}
\]
\[
\begin{split}
        & \leq Ch^{\beta} \int _0^t (t-\tau)^{-\beta/2 }(1+\|u(\tau)\|_{L_2(\Omega, H)})d\tau\\
        & + C\int_0^t \|u(\tau)-u_h(\tau)\|_{L_2(\Omega, H)}d\tau\\
         & \leq Ch^{\beta}+ C\int_0^t \|u(\tau)-u_h(\tau)\|_{L_2(\Omega, H)}d\tau.
     \end{split}
     \]
     Combining the estimates of $I_1(t) \text{ and }I_2(t)$, we obtain 
     \begin{equation}\label{errorv1}
     e_v(t)\leq Ch^{\beta }+C\int _0^t  \|u(\tau)-u_h(\tau)\|_{L_2(\Omega, H)}d\tau.
\end{equation}
     Above we have an estimate for $e_v(t)$.    
     We now find the error estimate between the exact solution $u(t)$ and the approximate solution $u_h(t)$. For $t \geq 0$, the error is given by 
     \[
     \begin{split}
         e_u(t)&=\|u(t)-u_h(t)\|_{L_2(\Omega, H)}\\
         &=\left\|(I+A)^{-1}v(t)-(\mathcal{P}_h+A_h)^{-1}v_h(t)\right\|_{L_2(\Omega, H)}\\
         &=\left\|(I+A)^{-1}v-(I+A)^{-1}v_h(t)+(I+A)^{-1}v_h(t)-(\mathcal{P}_h+A_h)^{-1}v_h(t)\right\|_{L_2(\Omega, H)}\\
         & \leq \|(I+A)^{-1}(v(t)-v_h(t))\|_{L_2(\Omega, H)}+ \left\|((I+A)^{-1}-(\mathcal{P}_h+A_h)^{-1})v_h(t)\right\|_{L_2(\Omega, H)}\\
         & \leq \|v(t)-v_h(t)\|_{L_2(\Omega, H)}+ \left\|((I+A)^{-1}-(\mathcal{P}_h+A_h)^{-1})v_h(t)\right\|_{L_2(\Omega, H)}\\
         &=e_v(t)+J(t),
     \end{split}
     \]
    where 
     \[J(t)=\left\|((I+A)^{-1}-(\mathcal{P}_h+A_h)^{-1})v_h(t)\right\|_{L_2(\Omega, H)}.
     \]
 To estimate $J(t)$, we use the connection between the resolvent of the generator and the corresponding $C_0$-semigroup   \cite[Chapter 3 Theorem 4.2]{apazy}.  For $t \geq 0$, we get
 \[
 \begin{split}
     J(t)&= \left\|((I+A)^{-1}-(\mathcal{P}_h+A_h)^{-1})v_h(t)\right\|_{L_2(\Omega, H)}\\
     & \leq \int _0^{\infty}e^{-\tau}\|(E_h(\tau)-E(\tau))v_h(t)\|_{L_2(\Omega, H)}d\tau\\
     &= \int _0^{\infty}e^{-\tau}\|(E_h(\tau)\mathcal{P}_h-E(\tau))v_h(t)\|_{L_2(\Omega, H)}d\tau\\
     &\leq Ch^{\beta}\int_0^{\infty} e^{-\tau}\tau^{-\beta/2} \|v_h(t)\|_{L_2(\Omega, H)}d\tau \\
     & \leq Ch^{\beta },
 \end{split}
 \]
 integral is finite  when $0 \leq \beta <2$.
 Here, we have used the identity $\mathcal{P}_hv_h=v_h$ for all $v_h \in V_h$ , bound of $v_h(t)$ ( proved in Lemma \ref{lemma1}) and the semigroup error estimate
 \[
 \|(E_h(t)\mathcal{P}_h-E(t))v\| \leq C h^{\beta}t^{-\beta/2}\|v\| \quad \text{ for } v\in H,\, t>0,
 \]
 which characterizes the strong convergence of the discrete semigroup $\{E_h(t)\}_{t \geq 0}$ to the continuous semigroup $\{E(t)\}_{t \geq 0}$. Combining the error estimate \eqref{errorv1} of $e_v(t)$ and $J(t)$, and using Gronwall's inequality , for $t \geq 0$  we get 
 \[
  \|u(t)-u_h(t)\|_{L_2(\Omega, H)}\leq Ch^{\beta}.
 \]
 This completes the proof of Theorem \ref{thm1}.
     \end{proof}
\begin{proof}[\bf Proof of Lemma \ref{lemma2}]
    We consider the solution of the stochastic parabolic component 
    \[
    v(t)=E(t)v_0+\int _0^t E(t-\tau)f((I+A)^{-1}v(\tau))d\tau+\int_0^t E(t-\tau)dW(\tau),\, \text{ for }\,t\geq 0.
    \]
    The solution of the stochastic parabolic component $v(t)$ is regular in time ( see  \cite[ Lemma 2.7]{Yan2} and \cite[ Proposition 3.4]{Jprintems}) i.e.
\[
\mathbb{E}\|v(t_2)-v(t_1)\|^2\leq C(t_2-t_1)^{\gamma},  \text{ for }\,  t_1\leq t_2\, \text{ and }\, 0\leq \gamma<\beta \leq 1.
\]
    From the change of variable $v(t)=(I+A)u(t) \text{ for } t \geq 0$, we have 
    \[
    \begin{split}
        \mathbb{E}\|u(t_2)-u(t_1)\|^2&=\mathbb{E} \|(I+A)^{-1}v(t_2)-(I+A)^{-1}v(t_1)\|^2\\
        &=\mathbb{E}\|(I+A)^{-1}(v(t_2)-v(t_1))\|^2\\
        &\leq\mathbb{E} \|v(t_2)-v(t_1)\|^2\leq C(t_2-t_1)^{\gamma}.
    \end{split}
    \]
 This completes the proof of Lemma \ref{lemma2}.   
\end{proof}
We now establish the mean-square stability of $V^n$, which plays a key role in the proof of Theorem \ref{thm2}.
\begin{proof}[\bf Proof of Proposition \ref{prop2}]
    From the equation \eqref{eqn2.9}, the solution $V^n$ admits the representation
    \[
    V^n=R(kA_h)^nV^0 +\sum_{j=0}^{n-1}\int_{t_j}^{j+1}R(kA_h)^{n-j}\mathcal{P}_hf(U^{j})d\tau +\sum_{j=0}^{n-1}\int_{t_j}^{t_{j+1}}R(kA_h)^{n-j}\mathcal{P}_hdW(\tau).
    \]
Applying the triangle inequality in $L_2(\Omega,H)$, we obtain
\[
\begin{split}
    \|V^n\|_{L_2(\Omega,H)} &\leq \|R(kA_h)^nV^0\|_{L_2(\Omega,H)} +\sum_{j=0}^{n-1}\int_{t_j}^{j+1}\|R(kA_h)^{n-j}\mathcal{P}_hf(U^{j})\|_{L_2(\Omega,H)}d\tau\\
    &+\| \sum_{j=0}^{n-1}\int_{t_j}^{t_{j+1}}R(kA_h)^{n-j}\mathcal{P}_hdW(\tau)\|_{L_2(\Omega,H)}\\
    &=:J_1+J_2+J_3.
\end{split}
\]
For the initial estimate, we have 
\[
J_1=\|R(kA_h)^nV^0\|_{L_2(\Omega,H)} =\|R(kA_h)^n\mathcal{P}_hv_0\|_{L_2(\Omega,H)} \leq C\|v_0\|_{L_2(\Omega, H)}.
\]
For the estimate of  $J_2$, we use the linear growth property of $f$ and Cauchy-Schwarz to get 
\[
\begin{split}
    J_2&=\sum_{j=0}^{n-1}\int_{t_j}^{j+1}\|R(kA_h)^{n-j}\mathcal{P}_hf(U^{j})\|_{L_2(\Omega,H)}d\tau\\
    &=\sum_{j=0}^{n-1}\int_{t_j}^{j+1}\left(\mathbb{E}\|R(kA_h)^{n-j}\mathcal{P}_hf(U^{j})\|^2\right)^{1/2}d\tau\\
    &\leq K_fC\sum_{j=0}^{n-1}\int_{t_j}^{j+1}\left(\mathbb{E}\|U^{j}\|^2\right)^{1/2}d\tau\\
            \end{split}
\]
\[
\begin{split}
    &=K_fCk\sum_{j=0}^{n-1}\left(\mathbb{E}\|U^{j}\|^2\right)^{1/2}\\
    & \leq K_fCk\sqrt{n}\left(\sum_{j=0}^{n-1}\mathbb{E}\|U^{j}\|^2\right)^{1/2}.
\end{split}
\]
Since $t_n=kn$ and $t_n \in [0,T]$, so we have
\[
J_2^2\leq C(T,f)k\sum_{j=0}^{n-1}\mathbb{E}\|U^{j}\|^2.
\]
For the estimate of $J_3$, we use It\^{o} isometry to obtain
\[
\begin{split}
    J_3^2&=\mathbb{E}\left\|\sum_{j=0}^{n-1}\int_{t_j}^{t_{j+1}}R(kA_h)^{n-j}\mathcal{P}_hdW(\tau)\right\|^2\\
    &=\sum_{j=0}^{n-1}\int_{t_j}^{t_{j+1}}\|R(kA_h)^{n-j}\mathcal{P}_hQ^{1/2}\|_{HS}d\tau\\
    &\leq Ckn Tr(Q)\leq CTr(Q).
\end{split}
\]
Combining the estimate of $J_1,\,J_2,\, \text{ and } J_3$, we get 
\[
\begin{split}
    \mathbb{E}\|V^n\|^2&\leq 3J_1^2+3J_2^2+3J_3^2\\
    & \leq 3\mathbb{E}\|v_0\|^2 +3C(T,f)k\sum_{j=0}^{n-1}\mathbb{E}\|U^j\|^2 + CTr(Q).
\end{split}
\]
From the full discretization of translation $U^j=(\mathcal{P}_h+A_h)^{-1}V^{j} \, \text{ for } j=0,1,\cdots,n$ and using discrete Gronwall's inequality, 
\[
\begin{split}
    \mathbb{E}\|V^n\|^2&\leq (3\mathbb{E}\|v_0\|^2 +3CTr(Q))e^{3(kn)}\\
    &\leq C(T,Q,v_0,f),\quad n=0,1,2,\cdots.
\end{split}
\]
This completes the proof of Proposition \ref{prop2}.
\end{proof}
To establish Theorem \ref{thm2}, we require certain regularity properties of the semigroup $\{E(t)\}_{t\geq 0}$. These are presented in the following lemma (see \cite{apazy,thomee}).
\begin{lemma}\label{lemma3}
 For any $ \nu \geq 0,\,0\leq\mu\leq 1$, there is a constant $C>0$ such that 
 \begin{align}
     &\|A^{\nu}E(t)\|\leq Ct^{-\nu}, \quad \text{ for }t>0, \label{estm1}\\
     &\|A^{\mu}(I-E(t))\|\leq Ct^{\mu}, \quad  \text{ for }t\geq 0\label{estm2}.
 \end{align}
\end{lemma}
Now we prove our second main result. First, we will compute the full discretization error between $V^n$ and $v(t_n)$ for all $ t_n \in [0, T]$. By using this error, we will find the full discretization error for the fourth-order stochastic pseudo-parabolic equation. 
\begin{proof}[\bf Proof of Theorem \ref{thm2}]
    From  the equation \eqref{eqn2.9}, we have 
    \[
    V^n=R(kA_h)^nV^0 +\sum_{j=0}^{n-1}\int_{t_j}^{j+1}R(kA_h)^{n-j}\mathcal{P}_hf(U^{j})d\tau +\sum_{j=0}^{n-1}\int_{t_j}^{t_{j+1}}R(kA_h)^{n-j}\mathcal{P}_hdW(\tau).
    \]
    From the stochastic parabolic component, we have
    \[
    v(t_n)=E(t_n)v_0+\int_0^{t_n}E(t_n-\tau)f(u(\tau))d\tau+\int_0^{t_n}E(t_n-\tau)dW(\tau).
    \]
    We define the error $e^n_v=V^n-v(t_n)$ and $F_n=R(kA_h)^n\mathcal{P}_h-E(t_n)$ to  get 
    \[
    \begin{split}
        e^n_v&=F_nv_0+\sum_{j=0}^{n-1}\int_{t_j}^{t_{j+1}}(R(kA_h)^{n-j}\mathcal{P}_hf(U^j)-E(t_n-\tau)f(u(\tau)))d\tau\\
        &+ \sum_{j=0}^{n-1}\int_{t_j}^{t_{j+1}}(R(kA_h)^{n-j}\mathcal{P}_h-E(t_n-\tau))dW(\tau)\\
        &=:Err_0+Err_d+Err_s,
    \end{split}
    \]
    where $Err_0$ is the error in the initial part, $Err_d$ is the error in the deterministic part, and $Err_s$ is the error in the stochastic part. We find the error estimate separately for each part.

    \textbf{Estimate for the initial error $Err_0$:} We use the Lemma $2.8$ of \cite{Yan2} to get 
    \[
\|Err_0\|_{L_{(\Omega,H)}}=\|F_nv_0\|_{L_2(\Omega,H)}\leq C(k^{\beta/2}+h^{\beta}),\quad \text{ for }0\leq \beta \leq 1.
    \]
       \textbf{Estimate for the deterministic error $Err_d$:}  For the deterministic error we have 
       \[
       \begin{split}(\mathbb{E}\|Err_d\|^2)^{1/2}&=\left(\mathbb{E}\left\|\sum_{j=0}^{n-1}\int_{t_j}^{t_{j+1}}(R(kA_h)^{n-j}\mathcal{P}_hf(U^j)-E(t_n-s)f(u(\tau)))d\tau\right\|^2\right)^{1/2}\\
       & \leq \sum_{j=0}^{n-1}\int_{t_j}^{t_{j+1}}\left(\mathbb{E}\|R(kA_h)^{n-j}\mathcal{P}_hf(U^j)-E(t_n-\tau)\|^2\right)^{1/2}d\tau\\
       &\leq \sum_{j=0}^{n-1}\int_{t_j}^{t_{j+1}}(\mathbb{E}\|R(kA_h)^{n-j}\mathcal{P}_h(f(U^j)-f(u(t_j)))\|^2)^{1/2}d\tau\\
       &+\sum_{j=0}^{n-1}\int_{t_j}^{t_{j+1}}(\mathbb{E}\|R(kA_h)^{n-j}\mathcal{P}_h(f(u(t_j))-f(u(\tau))\|^2)^{1/2}d\tau\\
       &+\sum_{j=0}^{n-1}\int_{t_j}^{t_{j+1}}(\mathbb{E}\|(R(kA_h)^{n-j}\mathcal{P}_h-E(t_n-t_j))f(u(\tau))\|^2)^{1/2}d\tau\\
       &+\sum_{j=0}^{n-1}\int_{t_j}^{t_{j+1}}(\mathbb{E}\|(E(t_n-t_j)-E(t_n-\tau))f(u(\tau))\|^2)^{1/2}d\tau\\
      &=:Err_{d,1}+Err_{d,2}+Err_{d,3}+Err_{d,4}.
       \end{split}
       \]
By using the global Lipschitz continuity property of $f$ to obtain
\[
\begin{split}
Err_{d,1}&=\sum_{j=0}^{n-1}\int_{t_j}^{t_{j+1}}(\mathbb{E}\|R(kA_h)^{n-j}\mathcal{P}_h(f(U^{j})-f(u(t_j)))\|^2)^{1/2}d\tau\\
&\leq Ck\sum_{j=0}^{n-1}(\mathbb{E}\|f(U^{j})-f(u(t_j))\|^2)^{1/2}\\
        \end{split}
\]
\[
\begin{split}
& \leq CK_fk\sum_{j=0}^{n-1}(\mathbb{E}\|U^j-u(t_j)\|^2)^{1/2}\\
&\leq C\sqrt{n}k\left(\sum_{j=0}^{n-1}\mathbb{E}\|U^j-u(t_j)\|^2\right)^{1/2}.
\end{split}
\]
Thus we have 
\[
Err_{d,1}^2\leq Ck\sum_{j=0}^{n-1}\mathbb{E}\|U^j-u(t_j)\|^2.
\]
For the estimate of $Err_{d,2}$, we use temporal continuity of $u$ ( Lemma \ref{lemma2}) to get 
\[
\begin{split}
Err_{d,2}&=\sum_{j=0}^{n-1}\int_{t_j}^{t_{j+1}}(\mathbb{E}\|R(kA_h)^{n-j}\mathcal{P}_h(f(u(t_j))-f(u(\tau)))\|^2)^{1/2}d\tau\\
&\leq CK_f \sum_{j=0}^{n-1}\int_{t_j}^{t_{j+1}}(\mathbb{E}\|u(t_j)-u(\tau)\|^2)^{1/2}d\tau\\
& \leq CK_f \sum_{j=0}^{n-1}\int_{t_j}^{t_{j+1}}(t_j-\tau)^{\gamma /2}d\tau\\
& \leq Ck^{\gamma /2}.
\end{split}
\]
For the estimate of $Err_{d,3}$, we apply Lemma $2.8$ of \cite{Yan2} and  linear growth condition \eqref{cond2} of  the function $f$ to obtain
\[
\begin{split}
    Err_{d,3}&=\sum_{j=0}^{n-1}\int_{t_j}^{t_{j+1}}(\mathbb{E}\|R(kA_h)^{n-j}\mathcal{P}_h-E(t_n-t_j)f(u(\tau))\|^2)^{1/2}d\tau\\
   &=\sum_{j=0}^{n-1} \int_{t_j}^{t_{j+1}}(\mathbb{E}\|F_{n-j}A^{(1-\beta)/2}A^{(\beta-1)/2}f(u(\tau))\|^2)^{1/2}d\tau\\
   &\leq \sum_{j=0}^{n-1}\|F_{n-j}A^{(1-\beta)/2}\|\int_{t_j}^{t_{j+1}}(\mathbb{E}\|A^{(\beta-1)/2}f(u(\tau))\|^2)^{1/2}d\tau\\
   &\leq Ck \sum_{j=0}^{n-1} \|F_{n-j}A^{(1-\beta)/2}\|\\
   & \leq C\sqrt{n}k\left(\sum_{j=0}^{n-1}\|F_jA^{(1-\beta)/2}\|^2\right)^{1/2}.
\end{split}
\]
From Lemma $2.8$ of \cite{Yan2}, we have 
\[
\begin{split}
Err_{d,3}^2&\leq C k\sum_{j=0}^{n-1}\|F_jA^{(1-\beta)/2}\|^2\\
&\leq Ck\sum_{j=0}^{n-1}\left(\sup_{v \neq 0}\frac{\|F_jA^{(1-\beta)/2}v\|}{\|v\|}\right)^2\\
        \end{split}
\]
\[
\begin{split}
& = C\sup_{v\neq 0}\frac{k\sum_{j=0}^{n-1}\|F_jA^{(1-\beta)/2}v\|^2}{\|v\|^2}\\
& \leq C \sup_{v \neq 0}\frac{(k^{\beta}+h^{2\beta})\|A^{(1-\beta)/2}v\|^2_{\beta-1}}{\|v\|^2}\\
& \leq C(k^{\beta}+h^{2\beta}).
\end{split}
\]
For the estimate of $Err_{d,4}$, we have
\[
\begin{split}
    Err_{d,4}&=\sum_{j=0}^{n-1}\int_{t_j}^{t_{j+1}}(\mathbb{E}\|(E(t_n-t_j)-E(t_n-\tau))A^{(1-\beta)/2}A^{(\beta-1)/2}f(u(\tau))\|^2)^{1/2}d\tau\\
    & \leq C\sum_{j=0}^{n-1}\int_{t_j}^{t_{j+1}}\|(E(t_n-t_j)-E(t_n-\tau)))A^{(1-\beta)/2}\|(\mathbb{E}\|A^{(\beta-1)/2}f(u(\tau))\|^2)^{1/2}d\tau\\
    &\leq C\sum_{j=0}^{n-1}\int_{t_j}^{t_{j+1}}\|(E(t_n-t_j)-E(t_n-\tau))A^{(1-\beta)/2}\|d\tau\\
    &\leq C\sqrt{n}\left(\sum_{j=0}^{n-1}\left(\int_{t_j}^{t_{j+1}}\|(E(t_n-t_j)-E(t_n-\tau))A^{(1-\beta)/2}\|d\tau\right)^2\right)^{1/2}\\
    & \leq C\sqrt{n}\sqrt{k}\left(\sum_{j=0}^{n-1}\int_{t_j}^{t_{j+1}}\|(E(t_n-t_j)-E(t_n-\tau))A^{(1-\beta)/2}\|^2d\tau\right)^{1/2}.
\end{split}
\]
Thus, by using Lemma \ref{lemma3}, we get 
\[
\begin{split}
    Err_{d,4}^2&\leq C\sum_{j=0}^{n-1}\int_{t_j}^{t_{j+1}}\|(E(t_n-t_j)-E(t_n-\tau))A^{(1-\beta)/2}\|^2d\tau\\
    &=  C\sum_{j=0}^{n-1}\int_{t_j}^{t_{j+1}}\|A^{1/2}E(t_n-t_j)A^{-\beta/2}(I-E(t_j-\tau))\|^2d\tau\\
    &\leq C \sum_{j=0}^{n-1}\int_{t_j}^{t_{j+1}}\|A^{1/2}E(t_n-t_j)\|^2\|A^{-\beta/2}(I-E(t_j-\tau))\|^2d\tau\\
   &\leq C\sum_{j=0}^{n-1}\|A^{1/2}E(t_n-t_j)\|^2\int_{t_j}^{t_{}j+1}(t_j-\tau)^{\beta}d\tau\\
   &\leq k^{\beta}C\sum_{j=0}^{n-1}(k\|A^{1/2}E(t_n-t_j)\|^2) \leq C k^{\beta}.
     \end{split}
\]
\textbf{Estimate for the stochastic error $Err_s$:} For the stochastic error, we have 
\[
\begin{split}
    \mathbb{E}\|Err_s\|^2&= \mathbb{E}\| \sum_{j=0}^{n-1}\int_{t_j}^{t_{j+1}}(R(kA_h)^{n-j}\mathcal{P}_h-E(t_n-\tau))dW(\tau)\|^2\\
    &=\sum_{j=0}^{n-1}\mathbb{E}\|\int_{t_j}^{t_{j+1}}(R(kA_h)^{n-j}\mathcal{P}_h-E(t_n-t_j)+E(t_n-t_j)-E(t_n-\tau))dW(\tau)\|^2\\
            \end{split}
\]
\[
\begin{split}
    & \leq \sum_{j=0}^{n-1}2 \int_{t_j}^{t_{j+1}}\|(R(kA_h)^{n-j}\mathcal{P}_h-E(t_n-t_j))Q^{1/2}\|_{HS}^2d\tau\\
    &+\sum_{j=0}^{n-1}2\int_{t_j}^{t_{j+1}}\|(E(t_n-t_j)-E(t_n-\tau))Q^{1/2}\|_{HS}^2d\tau\\
    &=2\sum_{j=0}^{n-1}\int_{t_j}^{t_{j+1}}\|(R(kA_h)^{n-j}\mathcal{P}_h-E(t_n-t_j))A^{(1-\beta)/2}A^{(\beta-1)/2}Q^{1/2}\|_{HS}^2d\tau\\
    &+2\sum_{j=0}^{n-1}\int_{t_j}^{t_{j+1}}\|(E(t_n-t_j)-E(t_n-\tau))A^{(1-\beta)/2}A^{(\beta-1)/2}Q^{1/2}\|_{HS}^2d\tau\\
    &\leq Ck\sum_{j=0}^{n-1}\|F_jA^{(1-\beta)/2}\|^2+ 2\sum_{j=0}^{n-1}\int_{t_j}^{t_{j+1}}\|(E(t_n-t_j)-E(t_n-\tau))A^{(1-\beta)/2}\|^2d\tau\\
    &\leq C(k^{\beta}+h^{2\beta})+Ck^{\beta} \leq C(k^{\beta}+h^{2\beta}).
\end{split}
\]
Combining all the estimates, $0 \leq \gamma< \beta \leq 1$ we have 
\begin{equation}\label{ferrorv}
\mathbb{E}\|e^n_v\|^2\leq C(k^{\gamma}+h^{2\beta})+Ck\sum_{J=0}^{n-1}\mathbb{E}\|U^j-u(t_j)\|^2.
\end{equation}
We now find the error estimate between the exact solution $u(t)$ and the fully discretized solution $U^n$. For any $t_n$ in $[0,T]$ the error is given by 
\[
\begin{split}
    \mathbb{E}\|U^n-u(t_n)\|^2&=\mathbb{E}\|(\mathcal{P}_h+A_h)^{-1}V^n-(I+A)^{-1}v(t_n)\|^2\\
    &=\mathbb{E}\|(\mathcal{P}_h+A_h)^{-1}V^n-(I+A)^{-1}v(t_n)\|^2\\
    &=\mathbb{E}\|(\mathcal{P}_h+A_h)^{-1}V^{n}-(I+A)^{-1}V^n+(I+A)^{-1}V^n-(I+A)^{-1}v(t_n)\|^2\\
    & \leq 2\mathbb{E}\|(\mathcal{P}_h+A_h)^{-1}V^n-(I+A)^{-1}V^n\|^2+\mathbb{E}\|V^n-v(t_n)\|^2.
\end{split}
\]
By using the connection between the resolvent of the generator and the corresponding $C_0$-semigroup \cite[ Chapter 3 Theorem 4.2]{apazy}, for $n\geq 0$, we get 
\[
\begin{split}
    \|((\mathcal{P}_h+A_h)^{-1}-(I+A)^{-1})V^{n}\|&\leq \int_{0}^{\infty}e^{-\tau}\|(E_h(\tau)-E(\tau))V^n\|d\tau\\
    &=\int_{0}^{\infty}e^{-\tau}\|(E_h(\tau)\mathcal{P}_h-E(\tau))V^n\|d\tau\\
    &\leq Ch^{\beta}\|V^n\|\int_{0}^{\infty}e^{-\tau}\tau^{-\beta/2}d\tau\\
    & \leq Ch^{\beta}\|V^n\|.
\end{split}
\]
By using the stability property of the fully discrete solution $V^n$, we have 
\[
\begin{split}
    \mathbb{E}\|((\mathcal{P}_h+A_h)^{-1}-(I+A)^{-1})V^{n}\|^2&\leq Ch^{2\beta}\mathbb{E}\|V^n\|^2\\
    &\leq Ch^{2\beta}.
\end{split}
\]
From the discrete Gronwall's inequality  and the estimate \eqref{ferrorv}    of $e^n_v$, we obtain 
\[
\begin{split}
    \mathbb{E}
    \|U^n-u(t_n)\|^2&\leq C(k^{\gamma}+h^{2\beta})+Ck\sum_{J=0}^{n-1}\mathbb{E}\|U^j-u(t_j)\|^2\\
    &\leq C(k^{\gamma}+h^{2\beta})e^{Ckn}\\
    & \leq  C(k^{\gamma}+h^{2\beta}).
\end{split}
\]
This completes the proof of Theorem \ref{thm2}.
\end{proof}
\section{Numerical Experiments}\label{S5}
We consider the following fourth-order stochastic pseudo-parabolic  equation 
\begin{equation}\label{eqn5.1}
    \begin{cases}
        du-d(\Delta u)-(  \Delta u-\Delta ^2u)dt=udt+dW \quad \text{ in } (0,1) \times (0,1),\\
        u(t,0)= u(t,1)=0 \hspace{4.8cm} t\in (0,1),\\
        \Delta u(t,0)=\Delta u(t,1)=0 \hspace{4cm} t\in (0,1),\\
        u(0,x)=\frac{1}{1+4\pi^2}\sin(2\pi x)\hspace{4.cm} x \in (0,1).
    \end{cases}
\end{equation}
From the change of variable $v=u-\Delta u$, the above equation converts \eqref{eqn5.1} into a stochastic coupled parabolic-elliptic system
\begin{equation}\label{eqn5.2}
    \begin{cases}
        dv-\Delta vdt=udt+dW \quad \text{ in }(0,1)\times (0,1),\\
        v=u-\Delta u \hspace{3.1cm} \text { in } (0,1)\times (0,1),\\
         v(t,0)=u(t,0)=0  \hspace{2cm} t\in (0,1),\\
         v(0,x)=\sin(2\pi x)\hspace{2.3cm} x\in (0,1).
          \end{cases}
\end{equation}
First, we consider the solution of the finite element approximation of the stochastic parabolic component of the stochastic coupled  parabolic-elliptic equation:
\[
v_h(t)=E_h(t)\mathcal{P}_hv_0+\int_0^tE_h(\tau)\mathcal{P}_hu_h(\tau)d\tau+\int_0^tE_h(\tau)\mathcal{P}_hdW(\tau) ,\quad t\in [0,1],
\]
where $u_h(t)$ is the solution of the elliptic component of the stochastic coupled parabolic-elliptic equation that is given by 
\[
u_h(t)=(\mathcal{P}_h+A_h)^{-1}v_h(t), \quad t\in [0,1].
\]
Consider a uniform partition of the time interval $[0,1]$, defined by  
$P_N=\{0=t_0 < t_1 < \cdots < t_N=1\}$,  
with time step size $k=1/N$. 
By using the semi-implicit time discretization and finite element approximation in the space with the mass matrix $M$ and diffusion matrix $K$, we get the discrete system for the stochastic parabolic component:
\[
MV^n+kA_hV^{n}=MV^{n-1}+k\mathcal{P}_h(M+K)^{-1}V^{n-1}+\mathcal{P}_h\Delta W^{n},\quad n\geq 1,
\]
where $\Delta W^n=W(t_n)-W(t_{n-1})$. The elliptic component is given by 
\[
U^{n}=(M+K)^{-1}V^{n-1},\quad n\geq1.
\]

\begin{figure}[htbp]
    \centering
    \begin{minipage}{0.45\textwidth}
        \centering
        \includegraphics[width=\linewidth]{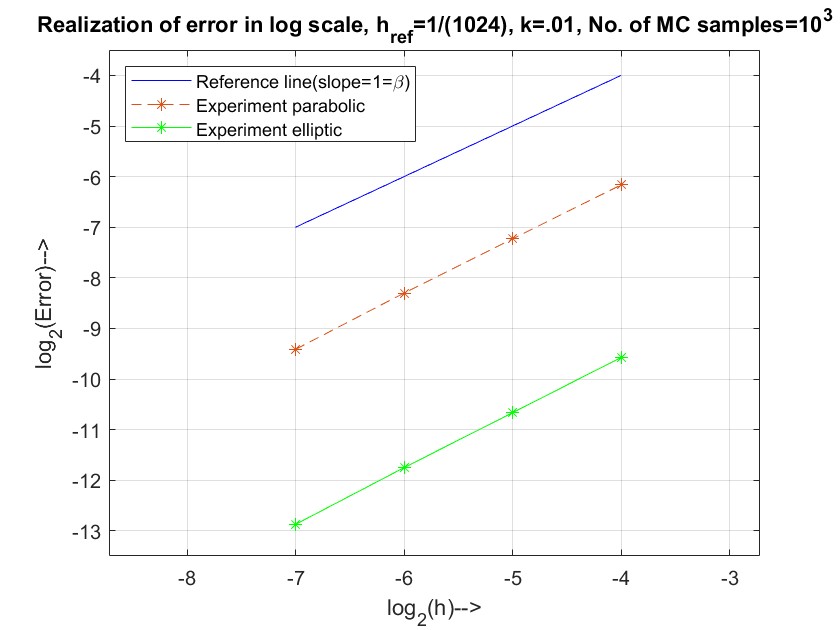} 
        \caption{The rate of strong convergence in $L_2$-norm with respect to  space discretized parameter}
        \label{fig1}
    \end{minipage}
    \hspace{0.05\textwidth} 
    \begin{minipage}{0.45\textwidth}
        \centering
        \includegraphics[width=\linewidth]{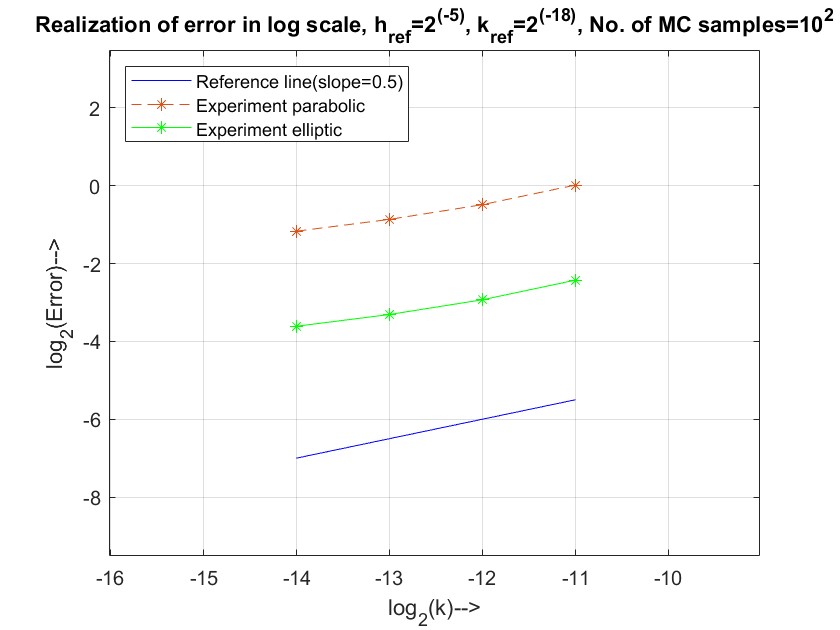} 
        \caption{The rate of strong convergence in $L_2$-norm with respect to  time discretized parameter}
        \label{fig2}
    \end{minipage}
\end{figure}

Using the full discretization scheme, we perform numerical computation for the example \eqref{eqn5.2} to support Theorem.  \ref{thm2}. Below, we present FIGURE \ref{fig1} and \ref{fig2}, which describe the rate of convergence for the equation \eqref{eqn5.2}.

Let $\{\lambda _i\}_{i=1}^{\infty}$ be eigen value of $A$ and set $Q=A^{-s}$. Then 
\[\|A^{(\beta-1)/2 }Q^{1/2}\|^2_{HS}=\|A^{(\beta-s-1)/2 }\|^2_{HS}=\sum_{i=1}^{\infty}\lambda_i^{\beta -s-1}\approx \sum_{i=1}^{\infty}i^{\frac{2}{d}(\beta -s-1)}
\]
which is finite if and only if 
$\beta < s+1-\frac{d}{2}$ where $d$ is the dimension of the domain $\mathcal{O}$. In the above example we have taken $d=1$. Hence, we require $\beta < s+\frac{1}{2}$.
In the numerical experiment, we have considered two cases:
\begin{enumerate}
    \item The rate of strong convergence of finite element approximations of fourth-order stochastic
pseudo-parabolic equation with respect to the space-discretized parameter $h:$ We take $\beta =1, \,\, d=1$, hence $s>1-\frac{1}{2}$. We choose $s=1/2+.0005$, $h_{ref}=1024$ (step length of space discretization for reference solution) and $k=.01$  k =0.01
(time
step to be fixed) and $10^3$ Monte-Carlo samples for sampling; see FIGURE \ref{fig1}.
\item The rate of strong convergence of time-discretized approximations of fourth-order stochastic pseudo-parabolic equations with respect to the time-discretized parameter $k$: We take $\gamma=0.98$, $d=1$, hence $s>0$. We choose $s=1/2+0.0005$, $h_{ref}=2^{-5}$ (step length to be
fixed ) and $k_{ref} = 2^{-17}$(step length of time discretization for reference solution) and
$10^2$ Monte-Carlo samples for sampling; FIGURE \ref{fig2}.

\end{enumerate}
\section{Conclusion}
In this work, we have analyzed the semi-discrete finite element approximation and the fully discrete scheme for the fourth-order stochastic pseudo-parabolic equation. 
Our results establish strong convergence rates for both spatial and temporal discretization parameters. 

An interesting direction for future research is the investigation of weak convergence for the numerical approximation of fourth-order stochastic pseudo-parabolic equations. 
Furthermore, we intend to investigate the strong convergence behavior of the scheme when the driving noise is a multiplicative $ Q$-Wiener process and L\'{e}vy process. 

\medskip
	\noindent
	
{\bf{\Large{Acknowledgement:}}} We express our gratitude to the Department of Mathematics and
	Statistics at the Indian Institute of Technology Kanpur for providing a conductive research
	environment. For this work,  M. Prasad is grateful
	for the support received through MHRD, Government of India (GATE fellowship). M. Biswas acknowledges support from IIT Kanpur. S. Bhar acknowledges the support from the SERB MATRICS grant (MTR/2021/000517), Government of India.
    
\vspace{1.5cm}
{\bf{\Large{Declarations:}}}\\
 \text{Ethics approval and consent to participate:} Not Applicable.\\
 \text{Consent for Publication: }Not Applicable.\\
 \text{Funding:} The authors declare that there was no funding for this work.\\
 \text{Data availability: } No data were used during the study.\\
\text{Author contributions: } Suprio Bhar has supervised the project and reviewed the manuscript. Mrinmay Biswas has selected the problem, reviewed the manuscript and contributed in the numerical simulation parts. Managala Prasad has done most of the technical part such as writing proofs and calculations, and has written the main manuscript.
	
	\bibliographystyle{plain}

\end{document}